\documentclass{article}
\usepackage{amsfonts,amssymb,amsbsy,amsmath,amscd,amstext,amsthm}
\usepackage{epsfig}

\theoremstyle{plain}
\newtheorem{lem}{Lemma}[section]
\newtheorem{thm}[lem]{Theorem}
\newtheorem{prop}[lem]{Proposition}
\newtheorem{coro}{Corollary}[lem]

\theoremstyle{definition}
\newtheorem{defn}[lem]{Definition}
\newtheorem{eg}{Example}

\theoremstyle{remark}
\newtheorem{rmk}[lem]{Remark}

\begin{document}

\title{A devil's staircase from rotations
        and irrationality measures for Liouville numbers}

\author{DoYong\ Kwon}
\date{}
\maketitle

\begin{abstract}
From Sturmian and Christoffel words we derive a strictly increasing
function $\Delta:[0,\infty)\rightarrow\mathbb{R}$. This function is
continuous at every irrational point, while at rational points,
left-continuous but not right-continuous. Moreover, it assumes
algebraic integers at rationals, and transcendental numbers at
irrationals. We also see that the differentiation of $\Delta$
distinguishes some irrationality measures of real numbers.
\end{abstract}

%
%
%
%
%
%
\section{Introduction}
For any real number $\alpha \in [0,1]$, we consider the dynamics of
rotation $$R_\alpha(\rho)= \rho+\alpha \mod 1.$$ Using an alphabet
$A=\{a,b\}$ with some integers $0\leq a<b$, the itinerary
$(R_\alpha^n(\rho))_{n\geq0}$ of some $\rho$ is recorded according
to a partition $$\mathcal{P}=\{[0,1-\alpha), [1-\alpha,1)\}\ \
\mathrm{or}\ \ \mathcal{P}'=\{(0,1-\alpha], (1-\alpha,1]\}.$$ In
other words, if $\rho\in [0,1]$ has traveled over $\mathcal{P}$,
then its infinite history $(s_{\alpha,\rho}(n))_{n\geq0}$ of
rotations is determined by the following rule:
$$s_{\alpha,\rho}(n)=
\left\{
  \begin{array}{ll}
    a, & \hbox{$R_\alpha^n(\rho)\in [0,1-\alpha)$;} \\
    b, & \hbox{$R_\alpha^n(\rho)\in [1-\alpha,1)$.}
  \end{array}
\right.
$$
In case that we adopt the partition $\mathcal{P}'$, a sequence
$(s'_{\alpha,\rho}(n))_{n\geq0}$ is also defined in a similar
fashion. One can observe that if $A=\{0,1\}$, then these infinite
words are obtained by the formulae
$$s_{\alpha,\rho}(n) = \lfloor \alpha(n+1) + \rho \rfloor -
\lfloor \alpha n +\rho \rfloor,\ \ s'_{\alpha,\rho}(n) = \lceil
\alpha(n+1) + \rho \rceil - \lceil \alpha n +\rho \rceil,$$ where
$\lfloor t \rfloor$ is the largest integer not greater than $t$, and
$\lceil t\rceil$ is the smallest integer not less than $t$. The
infinite words $s_{\alpha,\rho}$, $s'_{\alpha,\rho}$ are called
\textit{lower} and \textit{upper mechanical words} respectively with
\textit{slope} $\alpha$ and \textit{intercept} $\rho$. If the slope
$\alpha$ is irrational, then these infinite words are termed
\textit{Sturmian words}. For interesting properties of these words,
we refer the reader to \cite{Lo}.

While the rotations are `additive' dynamics, we now turn to
`multiplicative' dynamics. Let $\beta>1$ be a real number. The
$\beta$-\textit{transformation} $T_\beta:[0,1]\rightarrow [0,1)$
is a map defined by
$$T_\beta(x)=\beta x \mod 1.$$
Then the $T_\beta$-orbit of $x$ can be represented according to a
partition
$$\mathcal{P}=\{[0,1/\beta), [1/\beta,2/\beta),\ldots,
[(\lfloor\beta\rfloor-1)/\beta,\lfloor\beta\rfloor/\beta),
[\lfloor\beta\rfloor/\beta,1)\},$$ that is to say,
$$d_\beta (x):=(x_i)_{i\geq1}\
\mathrm{if\ and\ only\ if\ }x_i = \lfloor \beta T_\beta^{i-1} (x)
\rfloor.$$ We call $d_\beta(x)$ the $\beta$-\textit{expansion} of
$x$, and the $\beta$-\textit{shift} $S_\beta$ is the closure of $\{
d_\beta (x)|x\in[0,1) \}$ in the full shift. If $\beta$ is an
integer, then the $\beta$-expansion is nothing but the usual integer
base number representation and $S_\beta$ is the full shift over the
alphabet $A=\{0,1,\ldots,\beta-1\}$.

The $\beta$-expansion launched some study of real numbers. In terms
of the morphology of $d_\beta(1)$, Blanchard \cite{Bl} classified
real numbers greater than one into five classes $C_1$ to $C_5$. On
the other hand, Verger-Gaugry \cite{VG} focused on the patterns of
the consecutive zeros in $d_\beta(1)$, and defined classes
$Q_0^{(1)}$, $Q_0^{(2)}$, $Q_0^{(3)}$, $Q_1$ and $Q_2$. These new
points of view are quite different from the usual algebraic one. For
example, we do not know to which class $\frac{3}{2}$ belongs. See
\cite{Bl} and \cite{VG} for the precise definitions of each classes.

In \cite{CK}, irrational rotations were associated with
$\beta$-transformations in terms of their itinerary sequences, and
this connection naturally extends to rational rotations as explained
below. Chi and Kwon showed there that if $s$ is a Sturmian word of
irrational slope $\alpha\in (0,1)$ over an alphabet $A=\{a,b\}$ with
$0\leq a<b$, then there exists a unique $\beta\in(b,b+1)$ such that
$d_\beta(1)=E(s'_{\alpha,0})$ where $E$ is a morphism $0\mapsto a$,
$1\mapsto b$. Moreover they also proved that such $\beta$'s are
transcendental numbers, which provide continuum examples of
transcendental numbers in the class $C_3$. This ia a partial answer
to the question posed by Blanchard \cite{Bl}. Recently, a continuum
of transcendental numbers in the class $C_4$ were found in
\cite{AB3}.

From different contexts, Sturmian words also yielded transcendental
numbers in some numeration systems, e.g., in integer base expansions
\cite{FM} and in continued fraction expansions \cite{ADQZ}.
Recently, these two results were quite generalized by Adamczewski
and Bugeaud \cite{AB1,AB2}.

In the present paper, we will define a function
$\Delta:\mathbb{R}^{\geq0}\rightarrow \mathbb{R}$ by
$\Delta(\alpha)=\beta$ such that the $\beta$-expansion of $1$ is a
mechanical word of slope $\alpha$. Then the function enjoys some
devil's staircase-like properties:
\renewcommand{\theenumi}{\roman{enumi}}
\renewcommand{\labelenumi}{(\theenumi)}
\begin{enumerate}
  \item $\Delta$ is strictly increasing,
  \item at every irrational $\alpha>0$, $\Delta$ is continuous and
$\Delta(\alpha)$ is a transcendental number,
  \item at every rational $\alpha>0$, $\Delta$ is
left-continuous but not right-continuous and $\Delta(\alpha)$ is an
algebraic integer.
\end{enumerate}
For an irrational $\alpha$, one sees that $\Delta(\alpha)$ lies in
$C_3 \cap Q_0^{(1)}$. An algebraic study of $\Delta(\alpha)$ for
a rational $\alpha$ will appear in the subsequent paper \cite{Kw}.\\

We will also differentiate $\Delta$ wherever possible. Then apart
from the interest in its own right the differentiation of $\Delta$
distinguishes, in an unexpected manner, some number theoretic
properties of real numbers.

Let $t$ be a real number and $f$ be a positive increasing function
defined on positive integers. We consider a set
$$D_f =\left\{\frac{p}{q}\in \mathbb{Q}:
0<\left|t-\frac{p}{q}\right|<\frac{1}{f(q)} \right\}.$$ Over some
fixed class of functions, the supremum of functions $f$ that make
$D_f$ infinite is called the \textit{irrationality measure} of $t$.
One of classic and challenging questions in number theory is which
function $f$ allows $D_f$ to be finite or infinite. There is an
extensive literature seeking to find the irrationality measures of
some mathematical constants, e.g., \cite{Ge,Ro,Da,Ha90,Ha93,RV,Br}.
Recently, Adamczewski and Cassaigne \cite{AC} obtained an effective
upper bound of the irrationality measure of a real number whose
integer base representation is generated by a finite automaton.

The class of functions we usually consider is that of monomials or
that of exponential functions in $q$. The class of monomials defines
the \textit{irrationality exponent}, and the class of exponential
functions defines the \textit{irrationality base}. For instance, if
$t$ is a Liouville number then $D_f$ is an infinite set for any
monomial $f$ and the irrationality exponent of $t$ is infinite. In
Section \ref{Examples} we will also give examples of Liouville
numbers whose irrationality base is any real number greater than
$1$. By definition, the irrationality measure of $t$ measures how
well $t$ is approximable by rationals. In this sense Liouville
numbers are well approximable by rationals.

Let us divide, according to the differentiation of $\Delta$, the set
of positive irrational numbers into three sets by
$$\mathbb{R}^{\geq0}\setminus \mathbb{Q}
=\{x:\Delta'(x)=0\}\cup \{x:0<\Delta'(x)<\infty\}\cup \{x:\Delta\
\text{is\ not\ differentiable\ at\ }x \}.$$ We will see that this
trisection is related to the irrationality measure of $x$. Roughly
speaking, if $x\in \mathbb{R}^{\geq0}\setminus \mathbb{Q}$ is
`extremely' well approximable by rationals then $\Delta$ is not
differentiable at $x$. Otherwise, if $x$ is not `too' well
approximable by rationals then $\Delta'(x)$ exists and in most of
such cases it is equal to zero. Since almost every real number is
not well approximable by rationals \cite{Kh}, we can say that
$\Delta$ is flat almost everywhere or its almost all increases are
possible by discontinuous jumps only.

It is worthwhile to mention here that Bullett and Sentenac \cite{BS}
found a devil's staircase in a similar setting. The
$\beta$-expansions that are mechanical words were also considered.
But they fixed the base $\beta=2$, and instead investigated the
$\beta$-expansion $d_\beta(x)$ via varying $x$. More precisely, to
given $\alpha \in (0,1)$, they assigned the $T_\beta$-orbit whose
\textit{rotation number} is equal to $\alpha$. In fact they
considered, though the way of construction, its inverse function to
follow the usual definition of a devil's staircase: `the graph of a
non-constant function which is continuous, monotonic and locally
constant on a set of full measure.' See \cite{BS} for details. The
function $\Delta$ established in this paper is, however, nowhere
constant and thus does not fit the definition of a devil's
staircase. Its inverse function is indeed a devil's staircase in the
sense of \cite{BS}. But we will not take its inverse. This makes it
convenient to relate the differentiations of $\Delta$ with
irrationality measures.

%
%
%
%
%
%
\section{Christoffel words and a devil's staircase.}

We recall the notations and concepts on language theory. Lothaire's
book \cite{Lo} will be a complete alternative to this brief review.
Given a finite alphabet $A$, let $A\sp *$ (resp. $A\sp \mathbb{N}$)
be the set of finite (resp. infinite) words over $A$. A word $w \in
A\sp
* \cup A\sp \mathbb{N}$ is said to be a \textit{factor} (resp.
\textit{prefix}, \textit{suffix}) of a word $u\in A\sp * \cup A\sp
\mathbb{N}$ provided $u$ is expressed as $u=xwy$ (resp. $u=wy$,
$u=xw$) for some words $x$ and $y$. Note that if $u\in A\sp
\mathbb{N}$ then so is $y$. In case $A \subset \mathbb{N}$, the
usual \textit{lexicographic order} on $A\sp \mathbb{N}$ has a
natural extension to an order on $A\sp * \cup A\sp \mathbb{N}$ by
substituting any $x \in A\sp *$ with $x0^\omega:=x00\cdots$, even if
$0 \notin A$. For instance, if $x, y \in A\sp *$ and $z \in A\sp
\mathbb{N}$, then $x<y$ (resp. $y<z$, $z<y$) if and only if
$x0^\omega<y0^\omega$ (resp. $y0^\omega<z$, $z<y0^\omega$). This
total order is also called lexicographic order. A nonempty word
$u\in A\sp *$ is \textit{primitive} if $u=x^n$ for a nonempty $x$
implies $n=1$. For a word $u \in A\sp * \cup A\sp \mathbb{N}$, we
mean by $\mathrm{alph}(u)\subset A$ the set of letters appearing in
$u$. The set $A\sp\mathbb{N}$ is well endowed with a metric in a
sense that the metric generates the usual product topology of $A\sp
\mathbb{N}$ while $A$ has the discrete topology. For any $x, y \in
A\sp \mathbb{N}$, we define the distance between $x$ and $y$ by
$d(x,y)=2^{-n}$, where $n=\min \{k\geq 0| x_k \neq y_k \}$.

For the present we suppose $A=\{0,1\}$ and consider mechanical
words with their slope rational $\alpha=p/q \in[0,1]$,
$\gcd(p,q)=1$. It readily follows that $s_{\alpha,0}$ and
$s'_{\alpha,0}$ are purely periodic. To be more precise, we define
finite words $t_{p,q}$ and $t'_{p,q}$ by
$$t_{p,q}=a_0 \cdots a_{q-1},\quad t'_{p,q}=a'_0 \cdots
a'_{q-1},$$ where
$$a_i=\left\lfloor(i+1) \frac{p}{q} \right\rfloor-\left\lfloor i
\frac{p}{q} \right\rfloor,\ \ a'_i=\left\lceil(i+1)
\frac{p}{q}\right\rceil-\left\lceil i \frac{p}{q}\right\rceil.$$
Then one has $s_{\alpha,0}=t_{p,q}^\omega$ and
$s'_{\alpha,0}={t'_{p,q}}^\omega$. These words $t_{p,q}$,
$t'_{p,q}$ are said to be \textit{Christoffel words}. We observe
that $t_{1,1}=t'_{1,1}=1$, and that if $p/q \neq1$ then they are
factored as
$$t_{p,q}=0z_{p,q}1,\ \ t'_{p,q}=1z_{p,q}0,$$ for some word
$z_{p,q}$, called a \textit{central word}. It is easy to see that
all Christoffel words are primitive and that $z_{p,q}$ is a
\textit{palindrome word}, i.e., $z_{p,q}$ is equal to its reversal.
We recall that if $\alpha$ is irrational, then
$$s_{\alpha,0}=0c_\alpha\ \mathrm{and}\ s'_{\alpha,0}=1c_\alpha$$ for some
infinite word $c_\alpha$, called the \textit{characteristic word} of
slope $\alpha$. Before introducing our devil's staircase, we
investigate lexicographic order between some variants of Christoffel
words, which will play a crucial role in many contexts.

\begin{lem}
For $0<\rho,\rho' \leq1$, we have $s'_{\alpha,\rho}\leq
s'_{\alpha,\rho'}$ if $\rho< \rho'$.
\end{lem}

\begin{proof}
See \cite{Lo}.
\end{proof}

Note here that if the value $\alpha$ is restricted to irrationals
then we can say $s'_{\alpha,\rho}< s'_{\alpha,\rho'}$ if and only if
$\rho< \rho'$.

We denote by $\sigma$ the \textit{shift} of finite or infinite
sequences, and by $\{t\}$ the fractional part of $t$, i.e.,
$t=\lfloor t \rfloor + \{t\}$.

\begin{prop} \label{1z1}
For integers $p, q>0$, let $p/q \in(0,1)$ and $\gcd(p,q)=1$. Then
$1z_{p,q}1$ is lexicographically greater than all its proper
suffixes.
\end{prop}

\begin{proof}
Note $s'_{p/q,0}=(1z_{p,q}0)^\omega$ and put $\alpha=p/q$. We
define $\delta_0$ by
$$\delta_0=\min \{1-\{\alpha n\}: 1\leq n<q \}.$$
One finds that $0<\delta_0<1$ if $p$ and $q$ are relatively prime.
For $1\leq n<q$, let $\sigma^n (1z_{p,q}1)$ be a proper suffix of
$1z_{p,q}1$. Then this is also the prefix of $s'_{\alpha, \{\alpha
n\}+\delta_0/2}$. Since $\{\alpha n\}+\delta_0/2<1$, we conclude
from the lemma that
\begin{equation} \label{sfx1z1}
\sigma^n (1z_{p,q}1)< s'_{\alpha, \{\alpha n\}+\delta_0/2}\leq
s'_{\alpha,0}=s'_{\alpha,1}<1z_{p,q}1.
\end{equation}
\end{proof}

Noting that any proper suffix of
$1(z_{p,q}10)^\omega=1z_{p,q}1(0z_{p,q}1)^\omega$ can be expressed
as
$$\sigma^n (1(z_{p,q}10)^\omega)= s'_{p/q,
\{pn/q\}+\delta_0/2},\ n\geq1,$$ we know that Inequality
(\ref{sfx1z1}) shows the following.

\begin{coro} \label{1z10z1}
$1(z_{p,q}10)^\omega$ is lexicographically greater than all its
proper suffixes.
\end{coro}

In a $\beta$-shift $S_\beta$, Parry \cite{Pa} noticed that
$d_\beta(1)$ is quite distinguished in that $d_\beta(1)$ is greater
than all its proper suffixes. Moreover he also showed that this
property exhaustively determines whether a given sequence is
$d_\beta(1)$ for some $\beta>1$. Together with Parry's result,
Proposition \ref{1z1} guarantees the existence of a unique $\beta>1$
such that $d_\beta(1)=1z_{p,q}1$. This fact and \cite{CK} enable us
to define the next function.

\begin{thm}\label{DefOfDevil}
For integers $a$, $b$ with $0\leq a<b$, let $f_{a,b}(1)=b+1$. We
have a function $f_{a,b}:(0,1]\rightarrow [b,b+1]$ which satisfies
the following.
\renewcommand{\theenumi}{\alph{enumi}}
\renewcommand{\labelenumi}{{\rm(\theenumi)}}
\begin{enumerate}
  \item If $\alpha$ is irrational, then
       $d_{f_{a,b}(\alpha)}(1)=b c_\alpha$. In other words,
       $f_{a,b}(\alpha)$-expansion of $1$ is given by $b c_\alpha$.
  \item If $\alpha$ is equal to $p/q$ and $p$, $q$ are relatively prime,
       then $d_{f_{a,b}(\alpha)}(1)= b z_{p,q} b$.
\end{enumerate}
\end{thm}

\begin{rmk}
Note that the characteristic words and the central words involved in
Theorem \ref{DefOfDevil} are the same as when $A=\{0,1\}$ only
except that $0$ is replaced by $a$ and $1$ by $b$, that is,
$\mathrm{alph}(c_\alpha)=\mathrm{alph}(z_{p,q})=\{a,b\}$. Likewise,
we will read below $s_{\alpha,\rho}$ and $s'_{\alpha,\rho}$ as the
corresponding words over an alphabet $A=\{a,b\}$.
\end{rmk}

From now on we suppose $a=b-1$, so the alphabet is
$A=\{a,b\}=\{b-1,b\}$. We specify a function $\Delta:[0,\infty)
\rightarrow \mathbb{R}$ in such a way that for any positive integer
$b$, the restriction of $\Delta$ to the set $(b-1,b]$ is defined to
be the translation of $f_{b-1,b}$ by $b-1$, that is,
$$\Delta|_{(b-1,b]}(x) := f_{b-1,b}(x-b+1).$$ In addition we put
$\Delta(0):=1$. If we allowed the integer $a$ to be less than $b-1$,
then three different letters might appear in
$d_{f_{a,b}(\alpha)}(1)$, and hence many results to be obtained
below would fail.

\begin{lem} \label{limitofbeta}
Suppose that a sequence $(a_n)_{n\geq1}$ of real numbers is given to
satisfy $a_1\geq1$, $a_1 a_2 \cdots \neq 10^\omega$ and $0\leq
a_n\leq M$, $n=1, 2, \ldots$ for some $M$. Let $\zeta$ be the unique
positive root of a series equation
$$1=\sum_{n=1}^\infty a_n x^n.$$
We also assume that $\zeta_m$ is the unique positive root of
$$1=\sum_{n=1}^\infty b_{m,n} x^n,$$ where $0\leq b_{m,n}\leq M$,
$m,n=1, 2, \ldots$, and where $b_{m,i}=a_i$, $1\leq i\leq m$. Then
we have $$\lim_{m \to \infty}\zeta_m =\zeta.$$
\end{lem}
\begin{proof}
By equating two series, we get
$$a_1(\zeta_m-\zeta)+a_2(\zeta_m^2-\zeta^2)+\cdots+
a_m(\zeta_m^m-\zeta^m)= \sum_{n=m+1}^\infty a_n \zeta^n-
\sum_{n=m+1}^\infty b_{m,n} \zeta_m^n.$$ So one finds that
$$|\zeta_m-\zeta|\ |a_1+a_2(\zeta_m+\zeta)+\cdots+
a_m(\zeta_m^{m-1}+\cdots+\zeta^{m-1})|\leq
M\left(\sum_{n=m+1}^\infty \zeta^n + \sum_{n=m+1}^\infty \zeta_m^n
\right),$$ which is followed by
\begin{equation} \label{approx}
0\leq |\zeta_m-\zeta| \leq M \left(\frac{\zeta^{m+1}}{1-\zeta}
+\frac{\zeta_m^{m+1}}{1-\zeta_m}\right).
\end{equation}
The right-hand side of the inequality approaches zero as $m$ tends
to infinity because $0<\zeta_m<1-\varepsilon$ for some
$\varepsilon>0$.
\end{proof}

The next proposition determines the continuity of $\Delta$.

\begin{prop} For integers $a$, $b$ with $0\leq a=b-1$,
the function $f_{a,b}$ fulfills the following.
\renewcommand{\theenumi}{\alph{enumi}}
\renewcommand{\labelenumi}{(\theenumi)}
\begin{enumerate}
  \item $f_{a,b}$ is continuous at every irrational.
  \item $f_{a,b}$ is left-continuous but not right-continuous
       at every rational.
\end{enumerate}
\end{prop}

\begin{proof}
First we suppose that $\alpha_1$ is irrational, and
define $\delta_N$ by
\begin{equation} \label{sturmianpref}
\delta_N= \min \left\{ \frac{\{\alpha_1 n\}}{n}, \frac{1-\{\alpha_1
n\}}{n} :1\leq n\leq N \right\}.
\end{equation}
Since $\alpha_1$ is irrational, $\delta_N$ never vanishes for any
$N\geq1$. If $|\alpha-\alpha_1|<\delta_N$, then $s'_{\alpha,0}$ and
$s'_{\alpha_1,0}$ have a common prefix of length $N$. Hence Lemma
\ref{limitofbeta} shows that $|f_{a,b}(\alpha)-f_{a,b}(\alpha_1)|$
can be arbitrarily small if $N$ is sufficiently large.

For nonzero $\alpha_0=p/q\neq1$, $\gcd(p,q)=1$, we define $\delta_N$
by
\begin{equation} \label{lowchristpref}
\delta_N= \min \left\{ \frac{\{\alpha_0 n\}}{n}:
 1\leq n\leq N,\ n\not\equiv 0\mod q \right\}.
\end{equation}
Since $\gcd(p,q)=1$, one has $\delta_N>0$ for all $N\geq1$. If
$0\leq \alpha_0-\alpha <\delta_N$, then $s'_{\alpha,0}$ has a common
prefix of length $N$ with $s'_{\alpha_0,0}=(bz_{p,q}a)^\omega$.
Given $\varepsilon>0$, we can choose $N$ so that $0\leq
\alpha_0-\alpha <\delta_N$ implies
$$0\leq f_{a,b}(\alpha_0)-f_{a,b}(\alpha)<\varepsilon.$$ Here we use the fact
that if $\beta_0=f_{a,b}(\alpha_0)$ and $d_{\beta_0}(1)=\epsilon_0
\epsilon_1 \cdots \epsilon_{q-1}$ then an equation
$1=\sum_{n=0}^\infty {\epsilon'_n}/{\beta_0^{n+1}}$ also holds,
where $\epsilon'_n= \epsilon_n -1$ for $n\equiv -1 \mod q$ and
$\epsilon'_n= \epsilon_n$ otherwise.

To show that $f_{a,b}$ is not right-continuous, we assume
$\alpha>\alpha_0=p/q\neq1$. Then one finds
$$s'_{\alpha,0}>b (z_{p,q}ba)^\omega
>bz_{p,q}b>s'_{\alpha_0,0}.$$
With the aid of Corollary \ref{1z10z1}, we see that there exists
$\gamma>1$ such that
$$d_\gamma(1)= b (z_{p,q}ba)^\omega.$$ Then we
have $f_{a,b}(\alpha)> \gamma> f_{a,b}(\alpha_0)$.

It remains to prove left continuity at $1$. At first, one notes that
$f_{a,b}(1)=b+1$ is a positive root of $1=\sum_{n=1}^\infty b
x^{-n}$. So $0\leq 1-\alpha<1/N$ implies that $b^N$ is a prefix of
$s'_{\alpha,0}$. Therefore, $b+1-f_{a,b}(\alpha)$ tends to zero as
$N$ increases.
\end{proof}

The next theorem is easy consequences of the results above and
\cite{CK}.

\begin{thm}
We have the following.
\renewcommand{\theenumi}{\alph{enumi}}
\renewcommand{\labelenumi}{{\rm(\theenumi)}}
\begin{enumerate}
  \item $\Delta$ is strictly increasing.
  \item $\Delta(\mathbb{R^+})$ has Lebesgue measure zero.
  \item If $\alpha>0$ is irrational, then
       $$d_{\Delta(\alpha)}(1)=b c_\alpha,$$
       where $b=\lceil\alpha\rceil$ and $\mathrm{alph}(c_\alpha)=\{b-1,b\}$.
       In other words, $\Delta(\alpha)$-expansion of $1$
       is given by $b c_\alpha$. Moreover $\Delta(\alpha)$ is
       transcendental.
  \item If the fractional part $\{\alpha\}$ of $\alpha$ is equal to
       $p/q$ and $p$, $q$ are relatively prime, then
       $$d_{\Delta(\alpha)}(1)= b z_{p,q} b,$$
       where $b=\lceil\alpha\rceil$ and $\mathrm{alph}(z_{p,q})=\{b-1,b\}$.
       In this case, $\Delta(\alpha)$ is an algebraic
       integer.
  \item $\Delta$ is continuous at every irrational point.
  \item At every rational point, $\Delta$ is left-continuous but
       not right-continuous.
\end{enumerate}
\end{thm}

%
%
%
%
%
%
\section{Self-Christoffel numbers.}

In \cite{CK}, $\beta>1$ is called  a \textit{self-Sturmian number}
if $d_\beta(1)$ is a Sturmian word. Now we consider its counterpart
for Christoffel words.

\begin{defn}
For any rational $\alpha>0$, the value $\Delta(\alpha)$ is called a
\textit{lower self-Christoffel number}, and the right limit
$\Delta(\alpha+):=\lim_{x \to \alpha+} \Delta(x)$ called an
\textit{upper self-Christoffel number}.
\end{defn}

By the definition lower self-Christoffel numbers are all (simple)
beta-numbers and thus algebraic integers. This is also the case with
upper self-Christoffel numbers.

\begin{prop}
Suppose $\alpha=b-1+p/q$, $0< p\leq q$, $\beta_+ = \Delta(\alpha+)$
and $\gcd(p,q)=1$. Then
$$d_{\beta_+}(1)=b(z_{p,q}b(b-1))^\omega$$
if $p/q\neq1$, and
$$d_{\beta_+}(1)=(b+1)b^\omega$$
if $p=q=1$.
\end{prop}

\begin{proof}
If $p/q\neq1$, then we define $\delta_N$ by
\begin{equation} \label{upchristpref}
\delta_N= \min \left\{ \frac{1-\{p n/q\}}{n}:
 1\leq n\leq N,\ n\not\equiv 0\mod q \right\}.
\end{equation}
Since $p$ and $q$ are relatively prime, one notes that $\delta_N>0$
for all $N\geq1$. Now an inequality $0< \eta-\alpha <\delta_N$
implies that two words $s'_{\eta,0}$ and $b(z_{p,q}b(b-1))^\omega$
have a common prefix of length $N$. Now hence Corollary \ref{1z10z1}
and Lemma \ref{limitofbeta} prove
$d_{\beta_+}(1)=b(z_{p,q}b(b-1))^\omega$.

If $p=q=1$, then $0<\eta-b<1/N$ implies that
$\mathrm{alph}(s'_{\eta,0})=\{b+1,b\}$ and that $s'_{\eta,0}$ and
$(b+1)b^\omega$ have the identical prefix of length $N$. The same
reasoning as above also works.
\end{proof}

We now introduce canonical integer polynomials considered by Parry
\cite{Pa}, which have lower or upper self-Christoffel numbers as
zeros. Although the proof is immediate we state this fact as a
proposition. For a word $w=a_0 a_1 \cdots a_{n-1}$ with $a_i \in
\mathbb{Z}$, we mean by $\overrightarrow{w}$, a vector $(a_0,
\ldots, a_{n-1}) \in \mathbb{Z}^n$.

\begin{prop}
Suppose that $\alpha=b-1+p/q$, $0\leq p\leq q$, and $\gcd(p,q)=1$ if
$p\neq0$. If $\beta=\Delta(\alpha)$ and $\beta_+ = \Delta(\alpha+)$,
then $\beta$ is a root of an equation
$$x^q-\overrightarrow{b z_{p,q} b}\cdot (x^{q-1},x^{q-2},\ldots,1)=0,$$
and $\beta_+$ is a root of an equation
$$x^{q+1}-\overrightarrow{b z_{p,q} b}\cdot (x^q,x^{q-1},\ldots,x)-x+1=0.$$
\end{prop}

\begin{rmk}
Even if $q=1$, we will, for convenience, abuse $b z_{p,q}b$ for $b$
if $p=0$ and for $b+1$ if $p=1$. Bearing this convention in mind,
one can see that the preceding proposition covers the case where
$\alpha$ is an integer.
\end{rmk}

%
%
%
%
%
%
\section{Analysis on devil's staircase and irrationality measures}

Two polynomials in the previous section make it possible to get
some limit values of the devil's staircase $\Delta$.

\begin{lem} [\cite{Pa}] \label{atzero}
Let $\beta_n$ be the positive root of $1=x^{-1}+x^{-n}$. Then
\renewcommand{\theenumi}{\alph{enumi}}
\renewcommand{\labelenumi}{{\rm(\theenumi)}}
\begin{enumerate}
  \item $\lim_{n \to \infty} \beta_n =1$,
  \item $\lim_{n \to \infty} \frac{n}{\beta_n^n}=\infty$.
\end{enumerate}
\end{lem}

\noindent The above result is a key lemma for the next two
propositions.

\begin{prop}
The following hold.
\renewcommand{\theenumi}{\alph{enumi}}
\renewcommand{\labelenumi}{{\rm(\theenumi)}}
\begin{enumerate}
  \item $\lim_{\alpha \to 0+}\Delta(\alpha)=1$.
  \item $\Delta(b+)=\frac{b+2+\sqrt{b^2+4b}}{2}$
  for any integer $b\geq1$.
  \item $\lim_{\alpha \to
  \infty}(\Delta(\alpha)-\lceil\alpha\rceil)=1$.
\end{enumerate}
\end{prop}

\begin{proof}
Part (a) and (b) are immediate. To prove Part (c), we put
$$\beta_q(b):=\Delta(b-1+q^{-1}).$$
If $b=\lceil\alpha\rceil$ and $b-1+q^{-1}<\alpha$, then
$$1\geq \Delta(\alpha)-\lceil\alpha\rceil \geq \beta_q(b)-b.$$
Since $\beta_q(b)$ is a solution of $x^q-b
x^{q-1}=(b-1)(x^{q-2}+\cdots+1)+1$, one has
\begin{align*}
  \beta_q(b)-b &=(b-1)\left(\frac{1}{\beta_q(b)}+\cdots+\frac{1}{\beta_q(b)^{q-1}} \right)+
\frac{1}{\beta_q(b)^{q-1}}\\
   &=(b-1)\frac{1-\frac{1}{\beta_q(b)^{q-1}}}{\beta_q(b)-1}+\frac{1}{\beta_q(b)^{q-1}}.
\end{align*}
As $\lim_{b \to \infty}\beta_q(b)=\infty$ and
$b=\lfloor\beta_q(b)\rfloor$, we obtain the claim.
\end{proof}

\begin{figure}
\begin{center}
\epsfig{file=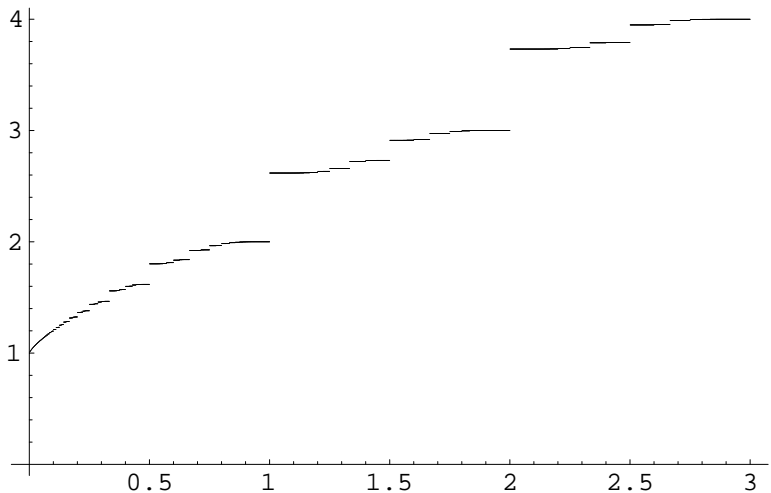} \caption{Devil's staircase}
\end{center}
\end{figure}

Although the graph of $\Delta$ is very ugly (see Figure 1), we can
do some calculus explicitly. This function has derivatives at almost
all irrational points, and left derivatives at all rational points.
Basically, $\Delta$ cannot have right derivative at rational points
since it is there discontinuous from the right, but we can say
something like right derivative as the following proposition says.
The proof employs a well-known fact, which is stated as a lemma for
future references.

\begin{lem}\label{LimitOfFunction}
$\lim_{x\to a} f(x)=l$ if and only if $\lim_{n\to \infty} f(a_n)=l$
for any sequence $(a_n)_{n\geq1}$ satisfying $a_n\neq a$ and
$\lim_{n\to\infty} a_n=a$.
\end{lem}

\begin{prop}\label{DerivativeAtRational}
The following hold for any rational $\alpha_0>0$.
\renewcommand{\theenumi}{\alph{enumi}}
\renewcommand{\labelenumi}{{\rm(\theenumi)}}
\begin{enumerate}
  \item $\lim_{\alpha \to 0+}\frac{\Delta(\alpha)-1}{\alpha}=\infty$.
  \item $\lim_{\alpha \to
  \alpha_0 -}\frac{\Delta(\alpha_0)-\Delta(\alpha)}{\alpha_0-\alpha}=0$.
  \item $\lim_{\alpha \to
  \alpha_0+}\frac{\Delta(\alpha)-\Delta(\alpha_0+)}{\alpha-\alpha_0}=0$.
\end{enumerate}
\end{prop}

\begin{proof}
(a) Let $q$ be a positive integer and $\beta_{2q}:=
\Delta((2^{-1}q^{-1})+)$ be the upper self-Christoffel number. For
any $\alpha \in (2^{-1}q^{-1}, q^{-1}]$, one finds
$$\frac{\Delta(\alpha)-1}{\alpha}>\frac{\beta_{2q}-1}{q^{-1}}.$$
Since $\beta_{2q}$ is the positive root of $x^{2q+1}-x^{2q}-2x+1=0$,
we have
$$\beta_{2q} -1=\frac{2\beta_{2q} -1}{\beta_{2q}^{2q}}.$$
Lemma \ref{atzero} now shows that
$$\lim_{q \to \infty}q(\beta_{2q}-1)
=\lim_{q \to \infty} \frac{q(2\beta_{2q} -1)}{\beta_{2q}^{2q}}
=\lim_{q \to \infty} \frac{q}{\beta_{2q}^{2q}}=\infty.$$ If a
positive sequence $(\alpha_n)_{n\geq1}$ satisfies
$\lim_{n\to\infty}\alpha_n=0$, then
$$\lim_{n\to\infty}\frac{\Delta(\alpha_n)-1}{\alpha_n}=\infty.$$ Since
$\bigcup_{q=1}^\infty (2^{-1}q^{-1}, q^{-1}]=(0,1]$, every
$\alpha_n$ lies in some interval $(2^{-1}q^{-1}, q^{-1}]$.\medbreak

Assume $\alpha_0=p/q\neq1$, $\gcd(p,q)=1$ and define $\delta_N$ as
Equation (\ref{lowchristpref}) for Part (b), and as Equation
(\ref{upchristpref}) for Part (c). In both cases, one observes that
$$\delta_N\geq \frac{1}{qN}\ \mathrm{and}\
\delta_{N+1}\geq \frac{1}{q(N+1)} \geq \frac{1}{2qN},$$ from which
it follows that the infinite union $\bigcup_{N=1}^\infty
((2qN)^{-1}, \delta_N]$ is an interval $(0,\{\alpha_0\}]$ or
$(0,1-\{\alpha_0\}]$ in each case.\\

\noindent (b) If $(2qN)^{-1}<\alpha_0-\alpha<\delta_N$, then
$s'_{\alpha,0}$ and $s'_{\alpha_0,0}$ have a common prefix of length
$N$. So Inequality (\ref{approx}) reads as follows:
\begin{align*}
\Delta(\alpha_0)-\Delta(\alpha) &\leq (b+1)^3
\left(\frac{1}{(\Delta(\alpha)-1)\Delta(\alpha)^N} +
\frac{1}{(\Delta(\alpha_0)-1)\Delta(\alpha_0)^N}\right)\\
&<  \frac{2(b+1)^3}{(c-1)c^N},
\end{align*}
for some $c>1$ independent of $N$. Thus one obtains
$$\frac{\Delta(\alpha_0)-\Delta(\alpha)}{\alpha_0-\alpha}<
\frac{4q N (b+1)^3}{(c-1)c^N}\rightarrow 0,\quad \mathrm{as\ }
N\rightarrow \infty.$$ Given a sequence $(\alpha_n)_{n\geq1}$
satisfying $\alpha_n<\alpha_0$ for $n\geq1$ and
$\lim_{n\to\infty}\alpha_n=\alpha_0$, we have now
$$\lim_{n\to\infty}
\frac{\Delta(\alpha_0)-\Delta(\alpha_n)}{\alpha_0-\alpha_n}=0.$$

\noindent (c) If $(2qN)^{-1}<\alpha-\alpha_0<\delta_N$, then
$s'_{\alpha,0}$ and $b(z_{p,q}b(b-1))^\omega$ have a common prefix
of length $N$. Now Inequality (\ref{approx}) implies that
\begin{align*}
\frac{\Delta(\alpha)-\Delta(\alpha_0+)}{\alpha-\alpha_0} &\leq 2qN
(b+1)^3 \left(\frac{1}{(\Delta(\alpha)-1)\Delta(\alpha)^N} +
\frac{1}{(\Delta(\alpha_0+)-1)\Delta(\alpha_0+)^N}\right)\\
&<  \frac{4qN (b+1)^3}{(\Delta(\alpha_0+)-1)\Delta(\alpha_0+)^N}
\rightarrow 0,\quad \mathrm{as\ } N\rightarrow \infty.
\end{align*}
If $\alpha_n>\alpha_0$ for $n\geq1$ and
$\lim_{n\to\infty}\alpha_n=\alpha_0$, then $$\lim_{n\to\infty}
\frac{\Delta(\alpha_n)-\Delta(\alpha_0)}{\alpha_n-\alpha_0}=0.$$ We
leave to the reader the proof when $\alpha_0$ is an integer.
\end{proof}

To find derivatives at irrational points is more involved, but uses
essentially the same technique as in finding one side derivatives at
rational points. We need some theory on Diophantine approximation.

We denote the regular continued fraction expansion of a real number
$t$ by
$$t=[a_0; a_1, a_2, \ldots]
   = a_0+\frac{1}{\displaystyle a_1
        +\frac{1}{\displaystyle a_2
        +\cdots}},$$
where the \textit{partial quotients} $a_i$ for $i\geq1$ are positive
integers. Then the $i$'th truncation $p_i/q_i
:=[a_0;a_1,\ldots,a_i]$ is called the $i$'th \textit{convergent} of
$t$.

Let $t$ be a real number and $p/q$ be a fraction in lowest terms
with $q>0$. Then we say that $p/q$ is:
\begin{itemize}
  \item a \textit{best approximation of the first kind} to the number $t$ if
  $$\left|t-\frac{p}{q} \right|<\left|t-\frac{a}{b} \right|$$
  holds for any fraction $a/b \neq p/q$ with $0<b \leq q$.
  \item a \textit{best approximation of the second kind} to the number $t$ if
  $$\left|q t-p \right|<\left|b t-a \right|$$
  holds for any fraction $a/b \neq p/q$ with $0<b \leq q$.
\end{itemize}
Because of an inequality
$$\left|t-\frac{p}{q} \right|
<\frac{b}{q} \left|t-\frac{a}{b}\right|\leq
\left|t-\frac{a}{b}\right|,$$ every best approximation of the second
kind is also a best approximation of the first kind. But the
converse is not true in general. Regular continued fractions and
Diophantine approximations have an intrinsic connection as the
following two propositions say. See \cite{Kh} for their proofs.

\begin{prop}\label{Lagrange}
Suppose that $t$ is a real number. Then
\renewcommand{\theenumi}{\alph{enumi}}
\renewcommand{\labelenumi}{(\theenumi)}
\begin{enumerate}
  \item every best approximation of the second kind to $t$ is a
  convergent of $t$,
  \item if $\{t\} \neq 1/2$, then every convergent of $t$ is
  a best approximation of the second kind to $t$.
\end{enumerate}
\end{prop}

The next proposition states that if $a/b$ is not a convergent of $t$
then it can not be too close to $t$.

\begin{prop}\label{convergents}
If $p/q$ is a fraction in lowest terms that fulfills the inequality
$$\left|t-\frac{p}{q} \right|<\frac{1}{2q^2},$$
then $p/q$ is a convergent of $t$.
\end{prop}

For a real number $\alpha$, let $\|\alpha\|$ be the smallest
distance from $\alpha$ to integers, namely,
$$\|\alpha\|:=\min_{n\in \mathbb{Z}}|\alpha-n|
=\min\{ \{\alpha\},1-\{\alpha\} \}.$$ If an irrational algebraic
$\alpha$ is of degree $n$, then there exists a constant
$c(\alpha)>0$ depending only upon $\alpha$ such that $\|
q\alpha\|>c(\alpha)q^{-n+1}$ for any nonzero $q\in \mathbb{Z}$. This
fact is known as Liouville's Theorem. Transcendental numbers passing
through this sieve are called \textit{Liouville numbers}, e.g.,
$\sum_{k=1}^{\infty}1/10^{k!}$. An irrational number $\alpha$ is
said to have \textit{irrationality exponent} $\mu(\alpha)$ if it is
given by
$$\mu(\alpha):=\sup\{\nu:\liminf_{q \to \infty} q^{\nu-1} \|q\alpha\|=0\}.$$
So we can say that Liouville numbers have their irrationality
exponent infinity and vice versa. It is well known that every
irrational number has its irrationality exponent greater than or
equal to $2$ and this value is in fact $2$ almost everywhere
\cite{Kh}. In particular, real numbers whose irrationality exponents
are equal to $2$ include all irrational algebraic numbers \cite{Ro}.
For our purpose we need to look at Liouville numbers in more detail.
For an irrational number $\alpha$, a number
$$\theta(\alpha)
:=\sup\{\lambda:\liminf_{q \to \infty} \frac{\lambda^q}{q}
\|q\alpha\|=0\}$$ is called the \textit{irrationality base} of
$\alpha$, which was coined by Sondow \cite{Son}. One can note the
following facts from the definitions.
\begin{itemize}
  \item \textit{$\mu(\alpha)<\infty$ implies $\theta(\alpha)=1$}.
  \item \textit{$\theta(\alpha)>1$ implies $\mu(\alpha)=\infty$}.
\end{itemize}
But there exists a real number $\alpha$ such that both
$\mu(\alpha)=\infty$ and $\theta(\alpha)=1$ hold at once. The
irrationality exponent and base can be reformulated as follows.

\begin{prop}\label{muandtheta}
For an irrational $\alpha$, suppose that $p_n/q_n$ is the $n$'th
convergent. Then
$$\mu(\alpha)=1-\liminf_{n \to \infty}\frac{\log\|q_n\alpha\|}{\log q_n},$$
$$\log\theta(\alpha)=-\liminf_{n \to \infty}\frac{\log\|q_n\alpha\|}{q_n}.$$
\end{prop}

\begin{proof}
Put $\mu=\mu(\alpha)$ and $\theta=\theta(\alpha)$. We know that
$\mu$ is the greatest upper bound of the set of all $\nu$ for which
$0<\|q\alpha\|<q^{-\nu+1}$ or
$$\nu<1-\frac{\log\|q\alpha\|}{\log q}$$
holds for infinitely many $q$. Thus one deduces that
$$\mu=\limsup_{q\to\infty}\left(1-\frac{\log\|q\alpha\|}{\log q}\right)
=1-\liminf_{q\to\infty}\frac{\log\|q\alpha\|}{\log q}.$$ Now we
claim that
$$\liminf_{q\to\infty}\frac{\log\|q\alpha\|}{\log q}
=\liminf_{n\to\infty}\frac{\log\|q_n\alpha\|}{\log q_n}.$$ If
$\tau=\liminf_{q\to\infty} \log\|q\alpha\|/\log q$, then given any
$\epsilon>0$,
$$\tau-\epsilon<\frac{\log\|q\alpha\|}{\log q}$$ holds for all sufficiently large $q$.
In particular,
$$\tau-\epsilon<\frac{\log\|q_n\alpha\|}{\log q_n}$$ holds for all
sufficiently large $n$. On the other hand, for infinitely many $q$
and $n$ , we have
$$\tau+\epsilon>\frac{\log\|q\alpha\|}{\log q} >\frac{\log\|q_n\alpha\|}{\log q_n},$$
where $q_{n-1}\leq q<q_n$. Hence we conclude
$\tau=\liminf_{n\to\infty} \log\|q_n\alpha\|/\log q_n$.

As for the irrationality base, $\theta$ is the greatest upper bound
of the set of all $\lambda$ for which $0<\|q\alpha\|<q/\lambda^q$ or
$$\log \lambda<\frac{\log q}{q}-\frac{\log\|q\alpha\|}{q}$$
holds for infinitely many $q$. The same reasoning as above proves
that
$$\log\theta =\limsup_{q\to \infty}\left(\frac{\log
q}{q}-\frac{\log\|q\alpha\|}{q}\right) = -\liminf_{q\to
\infty}\frac{\log\|q\alpha\|}{q}= -\liminf_{n\to
\infty}\frac{\log\|q_n\alpha\|}{q_n}.$$
\end{proof}

\begin{coro}\label{explicitmuandtheta}
$$\mu(\alpha)=1+\limsup_{n \to \infty}\frac{\log q_{n+1}}{\log q_n},$$
$$\log\theta(\alpha)=\limsup_{n \to \infty}\frac{\log q_{n+1}}{q_n}.$$
\end{coro}

\begin{proof}
These two equations follow from an inequality
$$\frac{1}{2q_{n+1}}<\|q_n\alpha\|<\frac{1}{q_{n+1}}.$$
\end{proof}

\begin{rmk}
Corollary \ref{explicitmuandtheta} was also noted by Sondow
\cite{Son}, but Proposition \ref{muandtheta} indeed verifies his
proof.
\end{rmk}

Although Liouville numbers are already well approximable by
rationals, we make the scale of well approximability more
minutely. We now trisect Liouville numbers according to their
approximability by rationals.

\begin{defn}
Let $\alpha$ be a Liouville number. We say that $\alpha$ is
\renewcommand{\theenumi}{\alph{enumi}}
\renewcommand{\labelenumi}{(\theenumi)}
\begin{enumerate}
  \item \textit{hypo-exponential} if
  $\mu(\alpha)=\infty$ and $\theta(\alpha)=1$,
  \item \textit{exponential} if
  $1<\theta(\alpha)<\infty$,
  \item \textit{hyper-exponential} if
  $\theta(\alpha)=\infty$.
\end{enumerate}
\end{defn}

In terms of $T_\beta$-orbit the value of $\Delta$ is, compared with
the other real numbers, quite exceptional. This point was
illuminated in \cite{CK} for irrationals, and in \cite{Kw} for
rationals.

\begin{lem}[\cite{CK,Kw}] \label{extremal}
Let $\alpha>0$ be real and $\beta=\Delta(\alpha)$. Then for all
integers $k\geq1$,
$$1-\frac1\beta< T_\beta^k(1)<1,$$
unless $T_\beta^k(1)=0$.
\end{lem}

If $d_\beta(1)$ is finite, say $d_\beta(1)=a_1 \cdots a_q$, then
$\beta$ is the unique positive root of an equation
$$1=\sum_{n=1}^q \frac{a_n}{x^n}.$$
But $\beta$ is also the root of a series equation
$$1=\sum_{n=1}^\infty \frac{a'_n}{x^n},\quad
a'_n=
\begin{cases}
    a_n -1 & \text{if\ $n\equiv 0\mod q$}, \\
    a_n & \text{otherwise}.
  \end{cases}
$$
For convenience' sake, we introduce a notation.
$$d_\beta(1-):= \lim_{x\to 1-} d_\beta(x)=
\begin{cases}
    d_\beta(1), & \text{if\ $d_\beta(1)$\ is\ infinite}, \\
    (a_1 \cdots a_{q-1}(a_q-1))^\omega, & \text{if\ $d_\beta(1)=a_1 \cdots a_q$}.
  \end{cases}
$$

The next lemma will show below that there exists some irrational
number where $\Delta$ is not differentiable.

\begin{lem}\label{lowerbound}
Given an irrational $\alpha>0$ and a real $\alpha_N>0$, let
$\beta=\Delta(\alpha)$ and $\beta_N=\Delta(\alpha_N)$. If
$d_\beta(1-)$ and $d_{\beta_N}(1-)$ have a common prefix of length
$N-1$ and their $N$'th letters are different, then there exist
constants $\delta\geq0$ and $c>0$ independent of $N$ such that the
following inequality holds:
$$|\beta-\beta_N|>\frac{c}{N (\beta+\delta)^N}.$$
\end{lem}

\begin{proof}
Let $d_\beta(1-)=a_1 a_2 \cdots$ and $d_{\beta_N}(1-)=b_1 b_2
\cdots$. Suppose at first $\alpha>\alpha_N$. We find then that
$$\frac{a_1}{\beta}+\cdots+\frac{a_{N-1}}{\beta^{N-1}}
+\frac{a_{N}}{\beta^{N}}+\frac{a_{N+1}}{\beta^{N+1}}+\cdots
=\frac{b_1}{\beta_N}+\cdots+\frac{b_{N-1}}{\beta_N^{N-1}}
+\frac{b_{N}}{\beta_N^{N}}+\frac{b_{N+1}}{\beta_N^{N+1}}+\cdots.$$
Since $a_1=b_1, \ldots, a_{N-1}=b_{N-1}$ and $a_N=b_N +1$, one has
\begin{align}\label{alpha>alphaN}
a_1 \left(\frac{1}{\beta_N}-\frac{1}{\beta} \right)&+ \cdots
+a_{N-1}\left(\frac{1}{\beta_N^{N-1}}-\frac{1}{\beta^{N-1}}\right)
+a_{N}\left(\frac{1}{\beta_N^{N}}-\frac{1}{\beta^{N}}\right) \\
     &> \frac{a_{N+1}}{\beta^{N+1}}+ \frac{a_{N+2}}{\beta^{N+2}}+\cdots
      >\left(1-\frac1\beta \right)\frac1{\beta^N},
\end{align}
where the last inequality follows from Lemma \ref{extremal}. On the
other hand, (\ref{alpha>alphaN}) is bounded above as
\begin{align*}
a_1 &\left(\frac{1}{\beta_N}-\frac{1}{\beta} \right) + \cdots
+a_{N}\left(\frac{1}{\beta_N^{N}}-\frac{1}{\beta^{N}}\right)\\
&=\left(\frac{1}{\beta_N}-\frac{1}{\beta} \right)
    \left(a_1+a_2\left(\frac{1}{\beta_N}+\frac{1}{\beta}\right)+\cdots+
       a_N\left(\frac{1}{\beta_N^{N-1}}+\cdots+\frac{1}{\beta^{N-1}}\right)
       \right)\\
&<bN\left(\frac{1}{\beta_N}-\frac{1}{\beta}\right)
\sum_{n=0}^\infty \frac{1}{\beta_N^i}=
\frac{bN(\beta-\beta_N)}{\beta(\beta_N-1)},
\end{align*}
where $b=\lceil \alpha \rceil$. Therefore we get
$$\beta-\beta_N>\frac{(\beta-1)(\beta_N-1)}{bN\beta^N}.$$

If $\alpha<\alpha_N$, then we have $b_N=a_N +1$. So one can show
\begin{align*}
\frac{bN(\beta_N-\beta)}{\beta_N(\beta-1)}&>
b_1\left(\frac{1}{\beta}-\frac{1}{\beta_N} \right)+ \cdots
+b_{N}\left(\frac{1}{\beta^{N}}-\frac{1}{\beta_N^{N}}\right) \\
     &> \frac{b_{N+1}}{\beta_N^{N+1}}+ \frac{b_{N+2}}{\beta_N^{N+2}}+\cdots
      >\left(1-\frac1{\beta_N}\right)\frac1{\beta_N^N}.
\end{align*}
Hence we have
$$\beta_N-\beta>\frac{(\beta-1)(\beta_N-1)}{bN\beta_N^N}
>\frac{(\beta-1)(\beta_N-1)}{bN(\beta+\delta)^N},$$
for some $\delta>0$.
\end{proof}

The next theorem enables us to differentiate $\Delta$ almost
everywhere. The proof depends on the completeness of real numbers.

\begin{thm} \label{DerivativeAtIrrational}
Suppose that $\alpha_0$ is an irrational number.
\renewcommand{\theenumi}{\alph{enumi}}
\renewcommand{\labelenumi}{(\theenumi)}
\begin{enumerate}
  \item If $\theta(\alpha_0)<\Delta(\alpha_0)$ then $\Delta'(\alpha_0)=0$.
  \item If $\theta(\alpha_0)>\Delta(\alpha_0)$ then $\Delta(x)$ is not differentiable at $x=\alpha_0$.
\end{enumerate}

\end{thm}

\begin{proof}
(a) For an integer $N\geq1$, we define $\delta_N$ by Equation
(\ref{sturmianpref}) or equivalently by
$$\delta_N= \min \left\{ \frac{\|\alpha_0 n\|}{n} :1\leq n\leq N
\right\}.$$ Let us fix a real $\lambda>1$ with
$\theta(\alpha_0)<\lambda<\Delta(\alpha_0)$. Then there exists an
integer $N_0\geq1$ such that $\delta_n>\lambda^{1-n}$ holds for
every $n\geq N_0$. Indeed, suppose otherwise that $L=\{n\in
\mathbb{N}: \delta_n\leq\lambda^{1-n}\}$ is an infinite set. Then so
is $\{ m_n:\delta_n=\|\alpha_0 m_n\|/m_n,\ m_n \leq n,\ n\in L\}$
--- otherwise, $\alpha_0$ is a rational number. For a real $\kappa$
with $\theta(\alpha_0)<\kappa<\lambda$ and $n\in L$, we have
$$\kappa^{m_n} \frac{\|\alpha_0 m_n\|}{m_n} =\kappa^{m_n}\delta_n
\leq \frac{\kappa^{m_n}}{\lambda^{n-1}}
=\frac{1}{\lambda^{n-m_n-1}}\frac{\kappa^{m_n}}{\lambda^{m_n}} \leq
\lambda\frac{\kappa^{m_n}}{\lambda^{m_n}} \rightarrow 0,\ \
\mathrm{as\ } m_n\rightarrow \infty,$$ which means
$\theta(\alpha_0)\geq \kappa$, a contradiction. Therefore if $n\geq
N_0$ then two intervals $(\lambda^{-n},\delta_n]$ and
$(\lambda^{-(n+1)},\delta_{n+1}]$ have a nonempty overlap, and so we
note $$\bigcup_{n=N_0}^\infty (\lambda^{-n},\delta_n]
=(0,\delta_{N_0}].$$ One sees that for all $\alpha$ satisfying
$\lambda^{-n}< |\alpha-\alpha_0|<\delta_n$, the words
$s'_{\alpha,0}$ and $s'_{\alpha_0,0}$ have a common prefix of length
$n$. Here we may assume $\lambda<\Delta(\alpha)$ for sufficiently
large $n$. Hence Lemma \ref{limitofbeta} guarantees some constant
$\gamma>1$ with $\lambda<\gamma<\Delta(\alpha_0)$, for which
\begin{align*}
|\Delta(\alpha)-\Delta(\alpha_0)| &\leq (b+1)^3
\left(\frac{1}{(\Delta(\alpha)-1)\Delta(\alpha)^n} +
\frac{1}{(\Delta(\alpha_0)-1)\Delta(\alpha_0)^n}\right)\\
&<  \frac{2(b+1)^3}{(\gamma-1)\gamma^n},
\end{align*}
where $b=\lceil \alpha_0\rceil$. Note that $\gamma$ is independent
of $n$. One thus finds that
\begin{equation}\label{DifferentiableAtHypoExponential}
\frac{|\Delta(\alpha)-\Delta(\alpha_0)|}{|\alpha-\alpha_0|}
<\frac{2(b+1)^3 \lambda^n}{(\gamma-1)\gamma^n} \rightarrow 0,\quad
\mathrm{as\ } n\rightarrow \infty.
\end{equation}
If $(\alpha_k)_{k\geq1}$ is a sequence such that $\alpha_k \neq
\alpha_0$ for $k\geq1$ and $\lim_{k\to\infty}\alpha_k=\alpha_0$ then
the value $|\alpha_k-\alpha_0|$ eventually lies in an interval out
of $\{(\lambda^{-n},\delta_n] : n\geq N_0 \}$, whence we get
$$\lim_{k\to\infty}
\frac{|\Delta(\alpha_k)-\Delta(\alpha_0)|}{|\alpha_k-\alpha_0|}=0.$$
Now Lemma \ref{LimitOfFunction} proves the claim.\medbreak

\noindent (b) Adopting $\delta_N$ as Part (a), we know that if
$\Delta(\alpha_0)<\lambda<\theta(\alpha_0)$ then the inequality
$\|\alpha_0 N\|/N<\lambda^{-N}$ and thus $\delta_N<\lambda^{-N}$
hold for infinitely many $N$. Let us fix $\lambda$ in the interval
$(\Delta(\alpha_0),\theta(\alpha_0))$ and collect such $N$'s by
$$Q=\{n\in\mathbb{N}:\delta_n \leq \frac{\|\alpha_0 n\|}{n}<\lambda^{-n}\}.$$
For sufficiently large
$q\in Q$, one can find an integer $p$ for which the inequality
$$\frac{\|\alpha_0 q\|}{q}=\left|\alpha_0-\frac{p}{q} \right|
<\frac{1}{\lambda^q}<\frac{1}{2q^2}$$ holds. Hence Proposition
\ref{convergents} shows that if $p$, $q$ are relatively prime then
$p/q$ is a convergent of $\alpha_0$. Otherwise if $d=\gcd(p,q)$ and
$p=dp', q=dq'$ then one has
$$\left|\alpha_0-\frac{p'}{q'} \right|
<\frac{1}{\lambda^{dq'}}<\frac{1}{\lambda^{q'}},$$ which implies
$q'\in Q$. We claim that infinitely many $q\in Q$ are indeed the
denominators of convergents of $\alpha_0$. If this is not the case,
then there exist an integer $q\in\mathbb{N}$ and an integer sequence
$(d_n)_{n\geq1}$ such that $\lim_{n \to \infty}d_n=\infty$ and $d_n
q\in Q$ for all $n\geq1$. In other words, we have for some integer
$p$,
$$\left| \alpha_0-\frac{p}{q}\right|<\frac{1}{\lambda^{d_n q}},\
n=1,2,\ldots.$$ But this inequality forces $\alpha_0$ to be equal to
$p/q$. This is a contradiction.

Now suppose that $q_1<q_2<\cdots$ are the denominators of the
convergents $p_i/q_i$ of $\alpha_0$ and that every $q_i$ lies in
$Q$. Then Proposition \ref{Lagrange} leads us to deduce that the
fraction $p_i/q_i$ is a best approximation of the second kind, and
hence of the first kind, i.e.,
$$\delta_{q_i}=\frac{\|\alpha_0 q_i\|}{q_i}\ \ \mathrm{and}\ \
\delta_{q_i}<\delta_{q_i-1}.$$ For each $i\geq1$, we define
$\alpha_{q_i}$ by
$$\alpha_{q_i}:=
  \begin{cases}
    \alpha_0-\delta_{q_i} & \text{if $\|q_i \alpha_0\|=q_i\alpha_0-p_i$}, \\
    \alpha_0+\min\{\frac{\delta_{q_i}+\delta_{q_i-1}}{2},2\delta_{q_i}\} &
\text{if $\|q_i \alpha_0\|=p_i-q_i\alpha_0$}.
  \end{cases}$$
Then one notes that
$$\lim_{i \to\infty}\alpha_{q_i}=\alpha_0\
\ \mathrm{and}\ \ |\alpha_{q_i}-\alpha_0| \leq
2\delta_{q_i}<\frac{2}{\lambda^{q_i}},$$ since $q_i$ is an element
of $Q$. Let us assume $\beta_0:=\Delta(\alpha_0)$,
$\beta_{q_i}:=\Delta(\alpha_{q_i})$ and $d_{\beta_0}(1)=a_1 \cdots
a_{q_i-1}a_{q_i}\cdots$. On the one hand if $\|q_i
\alpha_0\|=q_i\alpha_0-p_i$ then $\beta_{q_i}$ is a lower
self-Christoffel number, more precisely,
$$d_{\beta_{q_i}}(1-)=(a_1 \cdots a_{q_i-1}(a_{q_i}-1))^\omega.$$
On the other hand if $\|q_i \alpha_0\|=p_i-q_i\alpha_0$, then a word
$a_1 \cdots a_{q_i-1}(a_{q_i}+1)$ is a prefix of
$d_{\beta_{q_i}}(1-)$. In both cases, Lemma \ref{lowerbound} shows
that if $i$ is sufficiently large then
$$|\Delta(\alpha_{q_i})-\Delta(\alpha_0)|
>\frac{c}{q_i \gamma^{q_i}},$$
for some constants $c>0$ and $\Delta(\alpha_0)\leq\gamma<\lambda$.
Since $\gamma$ is eventually independent of $i$, we conclude the
following:
$$\frac{|\Delta(\alpha_{q_i})-\Delta(\alpha_0)|}{|\alpha_{q_i}-\alpha_0|}
>\frac{c\lambda^{q_i}}{2q_i \gamma^{q_i}} \rightarrow \infty,\quad
\mathrm{as\ } i\rightarrow \infty.$$
\end{proof}

In the case of $\theta(\alpha_0)>\Delta(\alpha_0)$, we can not say
that $\Delta'(\alpha_0)=\infty$ as in Proposition
\ref{DerivativeAtRational} (a). The proof of Part (b) only shows the
existence of a sequence $(\alpha_n)_{n\geq1}$ converging to
$\alpha_0$ for which $\lim_{n\to\infty}
|\Delta(\alpha_{n})-\Delta(\alpha_0)|/|\alpha_{n}-\alpha_0| =\infty$
holds.

Theorem \ref{DerivativeAtIrrational} covers two trivial cases where
$\theta(\alpha_0)=1$ and $\theta(\alpha_0)=\infty$.

\begin{coro}
Suppose that $\alpha_0>0$ is irrational.
\renewcommand{\theenumi}{\alph{enumi}}
\renewcommand{\labelenumi}{(\theenumi)}
\begin{enumerate}
  \item If $\alpha_0$ is either a non-Liouville or a
  hypo-exponential Liouville number, then $\Delta'(\alpha_0)=0$.
  \item If $\alpha_0$ is a hyper-exponential Liouville number,
  then $\Delta(x)$ is not differentiable at $x=\alpha_0$.
\end{enumerate}
\end{coro}

A real valued function $f(x)$ is said to satisfy the
\textit{Lipschitz condition} of \textit{order} $\eta$ at $x=x_0$ if
$\eta$ is the supremum of all $\zeta$ for which
$$|f(x)-f(x_0)|\leq C|x-x_0|^\zeta,\ \ C: \mathrm{constant}$$ holds on some open interval
containing $x_0$. Since the function $\Delta(x)$ is
non-differentiable on a dense subset of $[0,\infty)$, its second
derivative is nonsense everywhere. But this `local' Lipschitz
condition allows us to measure how flat $\Delta(x)$ is wherever it
is differentiable. The next theorem says that $\Delta$ is totally
flat almost everywhere.

\begin{thm}
Assume that $\alpha_0$ is an irrational number.
\renewcommand{\theenumi}{\alph{enumi}}
\renewcommand{\labelenumi}{(\theenumi)}
\begin{enumerate}
  \item If $\alpha_0$ is either a non-Liouville or a
  hypo-exponential Liouville number, then $\Delta(x)$ satisfies
  the Lipschitz condition of infinite order at $x=\alpha_0$.
  \item If $1<\theta(\alpha_0)<\Delta(\alpha_0)$,
  then $\Delta(x)$ satisfies the Lipschitz condition
  of order $\log \Delta(\alpha_0)/\log \theta(\alpha_0)$ at $x=\alpha_0$.
\end{enumerate}
\end{thm}

\begin{proof}
Suppose $\theta(\alpha_0)<\Delta(\alpha_0)$. Let $\lambda$ and
$\gamma$ be chosen as in the proof of Theorem
\ref{DerivativeAtIrrational} (a). So we have
$\theta(\alpha_0)<\lambda<\gamma<\Delta(\alpha_0)$. From Inequality
(\ref{DifferentiableAtHypoExponential}) it follows that
$$\frac{|\Delta(\alpha)-\Delta(\alpha_0)|}{|\alpha-\alpha_0|^\zeta}
<\frac{2(b+1)^3 \lambda^{\zeta n}}{(\gamma-1)\gamma^n}.$$ Now
$\lim_{n\to\infty} \lambda^{\zeta n}/\gamma^n=0$ if and only if
$\zeta< \log \gamma/\log \lambda$. Letting $\lambda \searrow
\theta(\alpha_0)$ and $\gamma \nearrow \Delta(\alpha_0)$, we obtain
the claim.
\end{proof}

%
%
%
%
%
%
\section{Examples}\label{Examples}

In this section we give some Liouville numbers that belong to each
of three classes, which make us believe that the derivative of
$\Delta$ assumes zero at every non-pathological constant. Some of
these examples can be found in \cite{Son}.

\begin{eg}
$\Delta'(\sqrt{2})=\Delta'(e)=\Delta'(\pi)=\Delta'(\log 2)
=\Delta'(\zeta(3))=\Delta'(\Gamma(1/4))=\Delta'(\Gamma(1/3))
=\Delta'(2^{\sqrt{2}}) =\Delta'(\log 3/\log 2)=\Delta'(e^\pi)=0$,
where
$$\zeta(s)=\sum_{n=1}^\infty \frac{1}{n^s},\ \text{for\ } s>1,$$
and
$$\Gamma(s)=\int_0^\infty e^{-t} t^{s-1}dt,\ \text{for\ } s>0.$$
\end{eg}

\begin{proof}
Liouville Theorem tells us that all algebraic irrational numbers
have finite irrationality exponents. Euler explicitly found the
continued fraction of $e$ by
$$e=[2;\overline{1,2m,1}]_{m\geq1}:=[2;1,2,1,1,4,1,1,6,1,\ldots].$$
With Corollary \ref{explicitmuandtheta}, one can show
$\mu(e)<\infty$. Or more explicitly it is known \cite{Da} that if
$\epsilon>0$ and $p,q\in \mathbb{N}$ with $q>q_0(\epsilon)$ for
some constant $q_0$ then
$$\left|e-\frac{p}{q} \right|
>\left(\frac12-\epsilon \right) \frac{\log\log q}{q^2 \log q}.$$
In \cite{Ha93,Ha90}, Hata showed that $\pi$ and $\log2$ have finite
irrationality exponents, explicitly,
$$\mu(\pi)\leq 8.0161, \ \mu(\log2)\leq3.8914.$$
As for $\zeta(3)$ the most up-to-date result due to Rhin and Viola
\cite{RV} is
$$\mu(\zeta(3))\leq5.513891.$$
In 2002, Bruiltet \cite{Br} obtained the explicit upper bounds of
irrationality measures for $\Gamma(1/4)$ and $\Gamma(1/3)$:
$$\mu(\Gamma(1/4))\leq 10^{143}, \ \mu(\Gamma(1/3))\leq 10^{151}.$$
Whether $e^\pi$ is a Liouville number is an open problem but the
known estimate given in \cite{Wa}
$$\left|e^\pi-\frac{p}{q} \right|>e^{-c(\log q)(\log\log q)},\
c=2^{60},\ q\geq3$$ is sufficient enough to show that the
irrationality base for $e^\pi$ is $1$. Indeed since for any
$\lambda>1$
$$\frac{\lambda^q}{e^{c(\log q)(\log\log q)}}
>\frac{\lambda^q}{e^{c(\log q)^2}}$$
tends to infinity as $q$ increases, we see that $\theta(e^\pi)=1$.
Gelfond Theorem \cite{Ge} implies that for any $\epsilon>0$ there
exist some constant $c_1,c_2>1$ such that
$$\left|2^{\sqrt{2}}-\frac{p}{q}\right|>(c_1 q)^{-(\log\log(c_1 q))^{5+\epsilon}}$$
and
$$\left|\frac{\log3}{\log2}-\frac{p}{q} \right|> e^{-(\log (c_2 q))^{4+\epsilon}}.$$
Along the same line as $e^\pi$, we can show that
$\theta(2^{\sqrt{2}})=\theta(\log3/\log2)=1$.
\end{proof}

A real number is said to be \textit{automatic} if its integer base
representation is generated by a finite automaton. For details the
reader is referred to \cite{AC}.

\begin{eg}
Let $\alpha>0$ be an automatic real number. Then
$\Delta'(\alpha)=0$.
\end{eg}

This example follows from the next theorem.

\begin{thm}[\cite{AC}]
Assume that
$$\{\alpha\}=\sum_{k=1}^{\infty}\frac{a_k}{b^k}$$ is an expansion
of the fractional part of $\alpha>0$ in base $b\in \mathbb{N}$, and
that the sequence $(a_k)_{k\geq1}$ is generated by a finite
automaton. Then $$\mu(\alpha)\leq C,$$ where the constant $C$
depends on the automaton giving $(a_k)_{k\geq1}$.
\end{thm}

\begin{eg}\label{FactorialSeries}
$\sum_{k=1}^{\infty}1/10^{k!}$ is a hypo-exponential Liouville
number.
\end{eg}

\begin{proof}
Let $\alpha_1=\sum_{k=1}^{\infty}1/10^{k!}$. It is well known that
$\alpha_1$ is a Liouville number. See for instance \cite{HW}. Let
$(p_n/q_n )_{n\geq1}$ be the sequence of the $n$'th convergents of
$\alpha_1$, and let
$$\frac{r_m}{s_m}:=\sum_{k=1}^{m}\frac{1}{10^{k!}}$$
be fractions in lowest terms, i.e., $\gcd(r_m,s_m )=1$. One can note
then that $s_m=10^{m!}$. From an inequality
\begin{align*}
  \alpha_1 &-\frac{r_m}{s_m} =\frac{1}{10^{(m+1)!}}+\frac{1}{10^{(m+2)!}}+\cdots \\
   & <\frac{1}{10^{(m+1)!}}\left(1+\frac{1}{10}+\frac{1}{10^2}+\cdots\right)
=\frac{10}{9\cdot 10^{(m+1)!}}<\frac{1}{2s_m^2},\ m = 2,3,4,\ldots,
\end{align*}
it follows that $(r_m/s_m)_{m\geq 2}$ is a subsequence of
$(p_n/q_n)_{n\geq1}$. We assume that a subsequence
$(q_{n_i})_{i\geq1}$ satisfies an equality
$$\lim_{i\to\infty}\frac{\log q_{n_i+1}}{q_{n_i}}=\log \theta(\alpha_1).$$
Then one can choose a (sub)sequence $\{s_{m_i} \}_{i\geq 1}$ of
$\{s_m \}_{m\geq 2}$ so that $s_{m_i}\leq q_{n_i}$ and
$q_{n_i+1}\leq s_{m_i+1}$ hold. Note here that we should allow some
$m_i$ to be equal to $m_{i+1}$. We have therefore
$$\log \theta(\alpha_1) =\lim_{i\to\infty}\frac{\log q_{n_i+1}}{q_{n_i}}
\leq\lim_{i\to\infty}\frac{\log s_{m_i+1}}{s_{m_i}}
=\lim_{i\to\infty}\frac{(m_i+1)!\log 10}{10^{m_i!}}=0.$$
\end{proof}

For a nonnegative integer $k$, we denote by $\mathcal{EXP}_k(x)$ the
$k$-fold power of $x$, for example,
$$\mathcal{EXP}_0(x)=1,\ \mathcal{EXP}_1(x)=x,\
\mathcal{EXP}_2(x)=x^x,\ \mathcal{EXP}_3(x)=x^{x^x}, \ldots.$$ Let
$\alpha=[a_0;a_1,a_2,\ldots]$ be irrational. From continued
fraction theory, we know that the more rapidly $a_n$ increases,
the more approximable by rationals $\alpha$ is. Although the
sequence $a_n=\mathcal{EXP}_n(10)$ increases extremely fast, it is
not enough for $\alpha_2=[a_0;a_1,a_2,\ldots]$ to be an
exponential or hyper-exponential Liouville number. Since
$q_{n+1}=a_{n+1}q_n+q_{n-1}$, we see that $a_{n+1}q_n< q_{n+1}<
(a_{n+1}+1)q_n$ and hence
$$\prod_{i=1}^n a_i <q_n <\prod_{i=1}^n(a_i+1).$$

\begin{eg}
$[\mathcal{EXP}_0(10);\mathcal{EXP}_1(10),\mathcal{EXP}_2(10),\ldots]$
is a hypo-exponential Liouville number.
\end{eg}

\begin{proof}
Suppose $a_n=\mathcal{EXP}_n(10)$ and
$\alpha_2=[a_0;a_1,a_2,\ldots]$. First $\alpha_2$ is actually a
Liouville number. This is because
\begin{align*}
\limsup_{n \to \infty} &\frac{\log q_{n+1}}{\log q_n}
 \geq\limsup_{n \to \infty}\frac{\log a_{n+1}+\log\prod_{i=1}^n a_i}{\log \prod_{i=1}^n(a_i+1)}
 \geq\limsup_{n \to \infty}\frac{\log a_{n+1}}{\log (2^n a_n^n)} \\
&=\limsup_{n \to \infty}\frac{a_n \log 10}{n(\log 2+ \log a_n)}
 \geq\limsup_{n \to \infty}\frac{10^{a_{n-1}}}{2n a_{n-1}\log10}=\infty.
\end{align*}
To show that $\alpha_2$ is hypo-exponential is also similar to the
above. We get
\begin{align*}
\limsup_{n \to \infty} &\frac{\log q_{n+1}}{q_n}
 \leq\limsup_{n \to \infty}\frac{\log \prod_{i=1}^{n+1}(a_i+1)}{\prod_{i=1}^n a_i}
 \leq\limsup_{n \to \infty}\frac{\log (2^{n+1}a_{n+1}^{n+1})}{\prod_{i=1}^n a_i} \\
 &=\limsup_{n \to \infty}\frac{(n+1)(\log2+a_n\log10)}{\prod_{i=1}^n
 a_i}=0.
\end{align*}
\end{proof}

From now on we construct exponential Liouville numbers.

Let a sequence $(a_n)_{n\geq1}$ be defined by
$$a_1=10,\ a_{n+1}=a_n !,\ n\geq1,$$
and we denote $a_{n+1}$ by ${10!}^n$. Then a number
$[0;10,10!,10!^2,10!^3,\ldots]$ still fails to be hyper-exponential.

\begin{eg}
$\alpha_3=[0;10,10!,10!^2,10!^3,\ldots]$ is an exponential Liouville
number. More precisely,
$$2.9221<\theta(\alpha_3)<7.9433.$$
\end{eg}

\begin{proof}
At first we note that
\begin{align*}
\log\theta(\alpha_3)
 &=\limsup_{n \to \infty}\frac{\log q_{n+1}}{q_n}
  =\limsup_{n \to \infty}\frac{\log (a_{n+1}q_n+q_{n-1})}{q_n}\\
 &=\limsup_{n \to \infty}\frac{\log a_{n+1} +\log q_n
   +\log\left(1+\frac{q_{n-1}}{a_{n+1}q_n}\right)}{q_n}\\
 &=\limsup_{n \to \infty}\frac{\log a_{n+1}}{q_n}
   =\limsup_{n \to \infty}\frac{\log (a_n!)}{q_n},
\end{align*}
and that
$$n\log n-n\leq \int_1^n \log x\ dx\leq \log n! \leq (n-1)\log n.$$
Hence we find
\begin{align*}
\log\theta(\alpha_3)
 &\leq \limsup_{n\to\infty} \frac{(a_n-1)\log a_n}{\prod_{i=1}^n a_i}
  \leq \limsup_{n\to\infty} \frac{(a_n-1)(a_{n-1}-1)\log a_{n-1}}{\prod_{i=1}^n a_i}\\
 & \leq \cdots \leq \limsup_{n\to\infty} \frac{\left(\prod_{i=1}^n (a_i-1)\right)\log
a_1}{\prod_{i=1}^n a_i}
 =\limsup_{n\to\infty} \log
 a_1\prod_{i=1}^n\left(1-\frac{1}{a_i}\right)\\
 &<\frac{9}{10}\log10.
\end{align*}
To the other direction, one can derive
\begin{align*}
\log\theta(\alpha_3)
 &\geq \limsup_{n\to\infty} \frac{a_n\log a_n-a_n}{\prod_{i=1}^n (a_i+1)}
  \geq \limsup_{n\to\infty}
    \frac{a_n a_{n-1}\log a_{n-1}-a_n a_{n-1}-a_n}{\prod_{i=1}^n (a_i+1)}\\
 &\geq \cdots \geq \limsup_{n\to\infty}
    \frac{(a_1\cdots a_n)\log a_1-[(a_1\cdots a_n)+(a_2\cdots a_n)+\cdots+a_n]}
         {\prod_{i=1}^n (a_i+1)}\\
 &=\limsup_{n\to\infty}
    \frac{\log a_1-[1+a_1^{-1}+(a_1 a_2)^{-1}+\cdots+(a_1\cdots a_{n-1})^{-1}]}
         {\prod_{i=1}^n (1+a_i^{-1})}.
\end{align*}
Now from
$$1+a_1^{-1}+(a_1 a_2)^{-1}+\cdots+(a_1\cdots a_{n-1})^{-1}
<\sum_{i=0}^\infty \frac{1}{10^i}=\frac{10}{9}$$ and from
$$\prod_{i=1}^n \left(1+\frac{1}{a_i}\right)
\leq \prod_{i=1}^n \left(1+\frac{1}{10^{2^{i-1}}}\right)
=\sum_{i=0}^\infty \frac{1}{10^i}=\frac{10}{9},$$ it follows that
$$\log\theta(\alpha_3) > \frac{9}{10}\log10 -1.$$
Gathering both results we have
$$10^{0.9} e^{-1}<\theta(\alpha_3)<10^{0.9}.$$
\end{proof}

\begin{rmk}
The author does not know what value the irrationality base of
$\alpha_3$ has.
\end{rmk}

Next, we give an exponential Liouville number that is represented by
series. Let
$$\alpha_4=\sum_{k=1}^\infty \frac{1}{\mathcal{EXP}_k(10)}
=\frac{1}{10}+\frac{1}{10^{10}}+\frac{1}{10^{10^{10}}}+\cdots.$$

\begin{eg}
$\theta(\alpha_4)=10$. So $\Delta'(9+\alpha_4)=0$ but $\Delta(x)$ is
not differentiable at $x=8+\alpha_4$.
\end{eg}

\begin{proof}
Let $(p_n/q_n)_{n\geq1}$ be the sequence of the $n$'th convergents
of $\alpha_4$, and let
$$\frac{r_m}{s_m}:=\sum_{k=1}^{m}\frac{1}{\mathcal{EXP}_k(10)}$$
be fractions in lowest terms. Then arguments similar to those of
Example 2 show that $s_m=\mathcal{EXP}_m(10)$ and
$$\alpha_4-\frac{r_m}{s_m} <\frac{10}{9\cdot \mathcal{EXP}_{m+1}(10)},$$
which in turn implies that $(r_m/s_m)_{m\geq 1}$ is a subsequence of
$(p_n/q_n)_{n\geq1}$. If $(q_{n_i})_{i\geq 1}$ satisfies
$$\lim_{i\to\infty}\frac{\log q_{n_i+1}}{q_{n_i}}=\log \theta(\alpha_4),$$
then as in Example \ref{FactorialSeries} there exists a
(sub)sequence $(s_{m_i})_{i\geq 1}$ such that $s_{m_i}\leq q_{n_i}$
and $q_{n_i+1}\leq s_{m_i+1}$ hold. So one can derive
$$\log \theta(\alpha_4) =\lim_{i\to\infty}\frac{\log q_{n_i+1}}{q_{n_i}}
\leq\lim_{i\to\infty}\frac{\log s_{m_i+1}}{s_{m_i}}
=\lim_{i\to\infty}\frac{s_{m_i}\log 10}{s_{m_i}}=\log10.$$ Since
$(s_m)_{m\geq 1}$ is a subsequence of $(q_n)_{n\geq1}$, Proposition
\ref{muandtheta} implies
$$\log\theta(\alpha_4)\geq -\liminf_{m \to \infty}\frac{\log\|s_m\alpha_4\|}{s_m}.$$
But we readily note
$$\|s_m\alpha_4\|
=\mathcal{EXP}_m(10)\sum_{k=m+1}^{\infty}\frac{1}{\mathcal{EXP}_k(10)}
\leq \frac{2\ \mathcal{EXP}_m(10)}{\mathcal{EXP}_{m+1}(10)}.$$
Thus we find
\begin{align*}
 \log\theta(\alpha_4)
  &\geq -\liminf_{m \to \infty}\frac{\log(2s_m/s_{m+1})}{s_m} \\
  &=-\liminf_{m \to \infty}\frac{\log2+\log s_m -\log s_{m+1}}{s_m}
   =\limsup_{m \to \infty}\frac{s_m\log10}{s_m}=\log10.
\end{align*}
\end{proof}

Let $\beta>1$ be fixed. We define two integer sequences
$(a_n)_{n\geq0}$ and $(q_n)_{n\geq-1}$ as follows. Firstly suppose
$q_{-1}=0$, $q_0=1$ and $a_0$ is any integer. Now define $a_n$ and
$q_n$ by
$$a_n=\left\lfloor \frac{\beta^{q_{n-1}}}{q_{n-1}}\right\rfloor,\
q_{n}=a_{n}q_{n-1}+q_{n-2}\ \ \text{for\ $n\geq1$}.$$ One can find
that these recursive equations uniquely determine $a_n$, $q_n$, and
also that $q_n$ is the denominator of the $n$'th convergent of
$[a_0;a_1,a_2,\ldots]$. The next example tells us that possible
irrationality bases exhaust all real numbers in $[1,\infty)$.

\begin{eg}
Let a sequence $(a_n)_{n\geq0}$ be as above and
$\alpha_5=[0;a_1,a_2,\ldots]$. Then $\theta(\alpha_5)=\beta$. So
$\Delta'(\lfloor\beta-1\rfloor+\alpha_5)=0$ but $\Delta(x)$ is not
differentiable at $x=\lfloor\beta-2\rfloor+\alpha_5$.
\end{eg}

\begin{proof}
By noting that
$$\beta^{q_n}-q_n\leq a_{n+1}q_n<q_{n+1}<(a_{n+1}+1)q_n\leq
\beta^{q_n}+q_n,$$ we have
$$\limsup_{n\to\infty}\frac{\log(\beta^{q_n}-q_n)}{q_n}
\leq \log\theta(\alpha_5)
\leq\limsup_{n\to\infty}\frac{\log(\beta^{q_n}+q_n)}{q_n},$$ that
is, $\theta(\alpha_5)=\beta$.
\end{proof}

The following are typical examples of hyper-exponential Liouville
numbers. We leave their proofs to the reader.

\begin{eg}
Suppose that $a_0$ is an arbitrary integer and for $n\geq1$ let
$$a_n=\mathcal{EXP}_n(n),\ b_n=\mathcal{EXP}_{2^n}(2^n).$$
Then $\alpha_6=[a_0;a_1,a_2,\ldots]$ and
$\alpha_7=\sum_{k=1}^\infty 1/b_n$ are hyper-exponential Liouville
numbers.
\end{eg}

We have seen some examples of hypo-exponential, exponential and
hyper-exponential Liouville numbers. In light of our problem, a
natural question arises whether the derivative of $\Delta$ can
assume a nonzero real number.\\

\noindent {\bf Question.} Does there exist a real $\alpha$ such that
$\theta(\alpha)=\Delta(\alpha)$? More specifically, is it possible
that $\Delta$ is differentiable at $\alpha$ and
$0<\Delta'(\alpha)<\infty$ ?

%
%
%
%
%
%

\noindent Department of Mathematics,\\
Yonsei University,\\
134 Shinchon-dong, Seodaemun-gu,\\
Seoul 120-749, Republic of Korea\\
E-mail: \textsf{doyong@yonsei.ac.kr}

\end{document}